\documentclass[11pt]{amsart}

\usepackage{amsmath,amsthm, amscd, amssymb, amsfonts}
\usepackage[all]{xy}

\usepackage[T2A, OT2, T1]{fontenc}

\usepackage[dvips, dvipsnames, usenames]{color}
\usepackage{graphicx}

\usepackage{t1enc}
\usepackage{fontsmpl}

\renewcommand{\_}[1]{_{( #1 )}}

\newcommand{\ce}[2]{C( #1, #2)}

\newcommand\ad{\operatorname{ad}}
\newcommand\sop{\operatorname{supp}}
\newcommand\Ad{\operatorname{Ad}}
\def\Aut{\operatorname{Aut}}

\newcommand\co{\operatorname{co}}
\def\End{\operatorname{End}}

\newcommand\Hom{\operatorname{Hom}}
\newcommand\Hopf{\operatorname{Hopf}}
\newcommand\Rep{\operatorname{Rep}}
\newcommand\Vect{\operatorname{Vec}}

\newcommand\id{\operatorname{id}}
\newcommand\Imm{\operatorname{Im}}
\newcommand\level{\operatorname{lv}}

\newcommand{\Hg}{\mathfrak H}
\newcommand{\Ag}{\mathfrak A}
\newcommand{\Bg}{\mathfrak B}
\newcommand{\Cg}{\mathfrak C}

\newcommand{\ku}{\Bbbk}

\newcommand{\N}{\mathbb N}

\def\C{\mathcal{C}}
\def\D{\mathcal{D}}
\def\M{\mathcal{M}}
\def\Ss{\mathcal{S}}

\newcommand{\leftexp}[2]{{\vphantom{#2}}^{#1}{#2}}

\newcommand{\ydh}{{}^{H}_{H}\mathcal{YD}}
\newcommand{\ydk}{{}^{K}_{K}\mathcal{YD}}

\numberwithin{equation}{section}\theoremstyle{plain}

\newtheorem{theorem}{Theorem}[section]

\newtheorem{lemma}[theorem]{Lemma}

\newtheorem{prop}[theorem]{Proposition}

\theoremstyle{definition}
\newtheorem{definition}[theorem]{Definition}
\newtheorem{example}[theorem]{Example}

\theoremstyle{remark}
\newtheorem{remark}[theorem]{Remark}

\newtheorem{problem}{Question}

\def\pf{\begin{proof}}
\def\epf{\end{proof}}

\begin{document}

\title[Examples of semisimple Hopf algebras]{Examples of extensions of Hopf algebras}

\author[Andruskiewitsch; M\"uller]
{Nicol\'as Andruskiewitsch, Monique M\"uller}

\address{FaMAF-CIEM (CONICET), Universidad Nacional de C\'ordoba,
Medina A\-llen\-de s/n, Ciudad Universitaria, 5000 C\' ordoba, Rep\'ublica Argentina.} \email{(andrus|mmuller)@famaf.unc.edu.ar}

\thanks{The work was partially supported by CONICET,
FONCyT-ANPCyT, Secyt (UNC)}

\begin{abstract}
We give some examples of, and raise some questions on, extensions of semisimple Hopf algebras. 
\end{abstract}

\maketitle


\section*{Introduction}

Semisimple Hopf algebras have been studied intensively in the last years because their representation categories 
have a rich structure-- they are fusion categories \cite{ENO1, ENO}. Conversely, one of the more fruitful approaches
to classification problems of semisimple Hopf algebras is through  fusion categories. However the basic notion
of extension of Hopf algebras is not categorical \cite{GN}, at least in a straightforward way. Here we explore this notion of extension.
First, we propose a definition of composition series of Hopf algebras.
After we shared this definition with her, S. Natale endorsed it with a suitable Jordan-H\"older theorem \cite{N3-JH}, see Theorem \ref{thm:jordan-holder}; 
cf. also \cite[Question 2.1]{bariloche}. Next we discuss briefly simple Hopf algebras and then present some new examples of length 2.
We raise questions along the text; some of them are very natural and perhaps not new, some of them are perhaps naive but we feel they might be useful.

The article is organized as follows. Section \ref{section:preliminaries} contains preliminaries. 
Section \ref{section:level-length} is devoted to  composition series and thoughts on simple Hopf algebras. 
In Section \ref{section:examples}, we discuss general ways to construct extensions, in particular Hopf algebras of length 2,
whose specific examples are in Section \ref{section:concrete examples}.

\smallbreak\subsubsection*{Acknowledgements} We thank C\'esar Galindo for pointing out to us Remark \ref{rem:galindo}, Alexei Davydov for 
pointing out to us Remark \ref{rem:davydov} and Sonia Natale for useful conversations and friendly remarks on a preliminar version.

\subsection*{Notations} 
Let $\ku$ be an algebraically closed field of characteristic 0. All  vector spaces, tensor products,   and algebras are over $\ku$. 
If $G$ is a group, then $\widehat{G} := \Hom (G, \ku^{\times})$. 
The notation $F\leqslant G$ means that $F$ is a subgroup of $G$, while $F \lhd  G$ means that $F\leqslant G$ is normal.
The standard basis of the group algebra    
$\ku G$ is $(g)_{g\in G}$, while $(\delta_g)_{g\in G}$ is the basis of the dual group algebra $\ku^G$ given by $\delta_g(h) = \delta_{g,h}$, $g, h\in G$. 
If $\triangleright$ is an action of a group $G$ on a set $F$, then we denote by $F^G$ the subset of $F$ of points fixed by $\triangleright$.

Let $H$ be a Hopf algebra, with coproduct $\Delta$, counit $\varepsilon$ and antipode $\Ss$. 
We use Heynemann-Sweedler's notation $\Delta(h) = h\_{1} \otimes h\_{2}$. The group of its group-like elements is $G(H)$
and the augmentation ideal is $H^+=\ker\varepsilon$.
If $\pi:H\rightarrow T$ is a Hopf algebra map, then its subalgebra of right coinvariants is 
$H^{\co\pi}=\{x\in H: (\id\otimes\pi)\Delta(x) = x\otimes 1\}$.  As usual, $\ydh$ denotes the category of Yetter-Drinfeld modules over $H$. An element $X=\sum_{i=1}^{r}X_i\otimes X^i\in A\otimes B$ is expressed $X=X_i\otimes X^i$ using the Einstein summation convention.

\section{Preliminaries}\label{section:preliminaries}

\subsection{Hopf algebra extensions} 

We write $K\leq H$ to express that $K$ is a Hopf subalgebra of $H$.
We denote a short exact sequence of Hopf algebras by 
$\xymatrix{ K\ar@{^{(}->}[r]^-\iota & H\ar@{->>}[r]^-\pi & T}$, meaning that $\iota$ and $\pi$ are Hopf algebra maps,
$\iota$ is injective, $\pi$ is surjective, $\ker \pi = H\iota(K^+)$ and $H^{\co\pi} = K$. We say that $H$ is an extension of $K$ by $T$.
From \cite[3.1]{Ma}, we get

\begin{remark}\label{rem:hopfsub}
If $H$ is an extension of $K$ by $T$ and $H' \leq H$, then $H'$ is an extension of $K' = K \cap H$ by $T' = \pi(H')$.
\end{remark}

Let $H$ a finite dimensional Hopf algebra. A Hopf subalgebra $K\leq H$ is  \textit{normal} if 
$\ad h(k) := h\_{1} k \Ss(h\_{2}) \in K$, for all $h\in H$, $k \in K$; we denote 
$K\lhd H$ and  $H/\hspace*{-.1cm}/ K:=H/HK^+$. If $K\lhd H$, then there is an exact 
sequence of Hopf algebras 
$\xymatrix{ K\ar@{^{(}->}[r]^-\iota & H\ar@{->>}[r]^-\pi & H/\hspace*{-.1cm}/ K }$, 
where $\iota$ is the inclusion and $\pi$ the canonical projection.
See \cite{AD, Ma, Sch}. 

\smallbreak
A finite dimensional Hopf algebra is  \textit{simple} if it does not contain proper normal Hopf subalgebras; then its
dual is simple too.  

\smallbreak
Let $K$, $T$ be Hopf algebras and $\rightharpoonup : T\otimes K \to K$, $\rho: T \to T\otimes K$, $\sigma: T\otimes T\to K$ 
and $\tau: T\to K \otimes K$ be linear maps. If $(\rightharpoonup, \rho, \sigma, \tau)$ is a \textit{compatible datum}, i.e., 
satisfies the conditions in \cite[\S 3]{A-canad} and $\sigma$, $\tau$ are invertible with respect to the convolution product, 
then there is a Hopf algebra structure $H=K{\,}\leftexp{\tau}{\#}_{\sigma}T$ on the vector space $K\otimes T$ given by 
$$(a\# g)(b\# h) = a(g\_1\rightharpoonup b)\sigma(g\_2,h\_1)\# g\_3h\_2,$$
$$\Delta(a\# g) = a\_1\tau(g\_1)_j\# \rho(g\_2)_i\otimes a_2\tau(g\_1)^j\rho(g\_2)^i\# g\_3.$$

An element $a\otimes g$ in $H$ is denoted by $a\# g$. Further, $H$ is an extension of $K$ by $T$; conversely, any cleft extension of $K$ by $T$ is $\simeq K{\,}\leftexp{\tau}{\#}_{\sigma} T$, 
for some compatible data $(\rightharpoonup, \rho, \sigma, \tau)$ \cite{AD}. All finite-dimensional extensions are cleft \cite{Sch}. The next is a particular case, that historically  appeared first.

\subsection{Abelian extensions}\label{subsection:abext} It is customary to call group algebras and their duals trivial Hopf algebras. 
The  first non-trivial examples of  finite-dimensio\-nal Hopf algebras are  the abelian extensions \cite{K, T, M}.
The input for the definition consists of

\begin{enumerate}\renewcommand{\theenumi}{\roman{enumi}}
\renewcommand{\labelenumi}{(\theenumi)}
	\item\label{item:matched-pair} A matched pair of finite groups $(F,G,\triangleleft,\triangleright)$;
	\item\label{item:cocycles} a pair of compatible cocycles $(\sigma, \tau)\in Z^2(F,(\ku^G)^\times) \times Z^2(G, (\ku^F)^\times)$,
	both normalized,	 denoted by
\begin{align*}
\sigma(x,y) &=\sum_{s\in G}\sigma_s(x,y)\delta_s,& x,y&\in F,& \tau(s,t)& =\sum_{x\in F}\sigma_x(s,t)\delta_x,& s,t&\in G,
\end{align*}
\end{enumerate}

Here \eqref{item:matched-pair} means that 
$\xymatrix{G & G\times F \ar[r]^-{\triangleright}\ar[l]_-{\triangleleft} & F}$  are (respectively, right and left) actions, such that
$s\triangleright xy = (s\triangleright x)((s\triangleleft x)\triangleright y)$ and   
$st\triangleleft x =(s\triangleleft(t\triangleright x))(t\triangleleft x)$, for all $s,t\in G$, $x,y\in F$.	
Given these data, there is a Hopf algebra structure
$H=\ku^{G}{\,}\leftexp{\tau}{\bowtie}_{\sigma}\ku F$ on the vector space $\ku^{G}\otimes\ku F$  
by $$(\delta_s\# x)(\delta_t\# y) = \delta_{s\triangleleft x, t}\sigma_s(x,y)\delta_s\# xy,$$
$$\Delta(\delta_s\# x) = \sum_{s=ab}\tau_x(a,b)\delta_a\# (b\triangleright x)\otimes \delta_b\# x.$$

Then $H$ fits into an extension  $\ku^G \hookrightarrow H\twoheadrightarrow \ku F$ 
and every extension of $\ku^G$ by $\ku F$ ari\-ses like this. 
Hopf algebras of this sort are named \emph{abelian extensions}.

\begin{lemma}\label{sub-abelian-is-abelian}
Hopf subalgebras of abelian extensions are abelian extensions.
\end{lemma}
\pf
Let $A$ fit into an exact sequence $\xymatrix@C-10pt{ \ku^G\ar@{^{(}->}[r]^-\iota & A\ar@{->>}[r]^-\pi & \ku F}$
and $B\leq A$. Since $\pi(B)$ is a Hopf subalgebra of $\ku F$,  $\pi(B)=\ku F'$, for some $F'\leqslant F$. 
Also, $B^{\co \pi|_{B}} = B\cap A^{\co\, \pi} = B\cap \ku^G$ is a Hopf subalgebra of $\ku^G$; thus there is a quotient group 
$G\rightarrow G'$ such that $B^{\co \pi|_{B}}=\ku^{G'}$. Then $B$ fits into the exact sequence 
$\xymatrix@C-10pt{\ku^{G'}\ar@{^{(}->}[r]& B\ar@{->>}[r] & \ku F'}$ by Remark \ref{rem:hopfsub}.
\epf

\subsection{Twisting}
There are  ways to obtain new Hopf algebras by altering the comultiplication or the multiplication.
The first appears in \cite{D} in the context of quasi-Hopf algebras; the second, dual to the first,
was first studied in \cite{DT}.

\subsubsection{Twisting the comultiplication}
Let $H$ be a Hopf algebra.
A twist for $H$  is $J\in H\otimes H$ invertible such that
\begin{align}\label{eq:twist}
\begin{aligned}
(\Delta\otimes\id)(J)(J\otimes 1)&=(\id\otimes\Delta)(J)(1\otimes J),
\\(\varepsilon\otimes\id)(J)&=(\id\otimes\varepsilon)(J)=1.
\end{aligned}
\end{align}

If $J$ is a twist, then the new Hopf algebra $H^J$ (called a twist of $H$) 
has underlying algebra $H^J=H$, comultiplication $\Delta^J(h)=J^{-1}\Delta(h)J$ 
and antipode $\Ss^J(h)=u^{-1}\Ss(h)u$, $h\in H$, where $u=m(\Ss\otimes\id)(J)$. 
If $J$ is a twist and $v\in H$ is invertible, then 
$J' =  \Delta(v)J(v^{-1} \otimes v^{-1})$ is again a twist, and $H^J \simeq H^{J'}$. 
One says that $J$ and $J'$ are  gauge-equivalent. 
Two Hopf algebras $H$ and $K$ are gauge-equivalent 
if and only if the tensor categories $\Rep H$ and $\Rep K$ are equivalent \cite{S}.

Twists in a group algebra $H =\ku N$, where $N$ is a finite group, are classified, up to gauge equivalence, 
by classes of pairs $(S,\omega)$, where $S\leqslant N$ is solvable, 
$\vert S\vert$ is a square and $\omega \in H^2(S,\ku^\times)$ is a non-degenerate 2-cocycle on $S$. See \cite{Mov, EG}.
Namely, if $J$ is a twist for $H$, then $S\leqslant N$ is the  subgroup  minimal for $J$, i.~e., 
the components of $J^{-1}_{21}J$ span $\ku S$; and $J$ determines $\omega$. 
For instance, if $S$ is abelian, 
then the twist corresponding to $(S, \omega)$ is given by $J=\sum_{\chi, \eta\in \widehat{S}}\omega(\chi,\eta)e_\chi\otimes e_\eta$, 
where $e_\chi=\dfrac{1}{|S|}\sum_{h\in S}\chi(h^{-1})h$.

\begin{theorem}\label{th_eg-triangular} \cite[Theorem 2.1]{EG2}
If $H$ is a semisimple triangular Hopf algebra, then it is isomorphic to a twist of a group algebra. \qed
\end{theorem}

We shall need the following fact.

\begin{lemma}{\cite[2.6]{GN}}\label{lemma:cocomm Twisted group algebra}
Let $J\in\ku N\otimes\ku N$ be the twist associated to the pair $(S,\omega)$, where $S$ is the minimal subgroup of $J$. 
Then $(\ku N)^J$ is cocommutative if and only if $S\lhd N$, $S$ is abelian and $\omega$ is \textrm{ad}$N$-invariant in $H^2(S,\ku^\times)$. \qed
\end{lemma}

\subsubsection{Twisting the multiplication}
Dually, a cocycle for a Hopf algebra $H$ is a linear map $\sigma: H\otimes H\to \ku$, invertible for the convolution product and satisfying 
the requirements dual to \eqref{eq:twist}. There is a new Hopf algebra $H_\sigma$ (called a  cocycle-twist of $H$)  with
underlying coalgebra $H$ and multiplication $m_\sigma = \sigma^{-1}*m * \sigma$. 
If $H$ is finite-dimensional and $\sigma$ is a cocycle for $H$, then $J = \sigma^*(1) \in H^* \otimes H^*$ is a twist and 
$(H_\sigma)^* \simeq (H^*)^J$. It follows at once from Theorem \ref{th_eg-triangular} that
\emph{if $H$ is a semisimple cotriangular Hopf algebra, then it is isomorphic to a cocycle-twist of a dual group algebra.}

\subsection{Weakly group-theoretical Hopf algebras}\label{subsec:wgt} 
Semisimple Hopf algebras can be studied through their representations; the category $\Rep H$ is a fusion category
\cite{ENO1}. We discuss briefly the interplay between classes of fusion categories and extensions.

Two fusion categories $\C$ and $\D$ are Morita-equivalent if there exists an indecomposable $\C$-module category $\M$
such that $\D$ is equivalent to $\End_{\C} (\M)$ \cite{ENO1}. 
This relation settles a basic reduction in the classification program of fusion categories 
and has the following counterpart: two semisimple Hopf algebras $H$ and $K$ are Morita-equivalent 
if  $\Rep H$ and $\Rep K$ are Morita-equivalent. 
Taking $\M = \Vect$, we see that $H$ and $H^*$ are Morita-equivalent.
Thus, gauge equivalence implies Morita equivalence, but not vice versa.

A fusion category is pointed if all simple objects are invertible; all pointed fusion categories are of the form
$\Vect^G_{\omega}$, that is, categories of $G$-graded vector spaces (for some finite group $G$) with associator
induced by $\omega \in H^3(G, \ku^{\times})$.
A fusion category Morita-equivalent to a pointed one is called \emph{group-theoreti\-cal} \cite{O}. 
A semisimple Hopf algebra $H$ is group-theoretical  when Rep $H$ is so; hence
the class of group-theoretical Hopf algebras is stable under twisting.
Abelian extensions are group-theoretical \cite{N}. 
It was conjectured that every semisimple Hopf algebra is group-theoretical \cite{ENO1}, but counterexamples were found in 
\cite{Nik2}: there is a non-group-theoretical Hopf algebra $H_p$ of dimension $4p^2$ for each odd prime $p$. 
Indeed, $H_p$ fits into an exact sequence  $\ku^{\mathbb{Z}/2}\hookrightarrow H_p \twoheadrightarrow (\ku G)^J$, 
where $G= (\mathbb{Z}/p\times\mathbb{Z}/p)\rtimes\mathbb{Z}/2$ 
and $J$ is a non-trivial twist in $\ku G$. Thus the class of group-theoretical Hopf algebras is not stable under extensions. 
More general examples of non-group-theoretical Hopf algebras were described in \cite{GNN}. 

\smallbreak
The notions of weakly group-theoretical fusion categories and Hopf algebras were introduced in \cite{ENO};
\cite[Question 2]{ENO} asks whether any semisimple Hopf algebra is weakly group-theoretical.

\begin{problem}\label{pbm:wgt-ext} \cite{N-pmu} Is  an extension 
of weakly group-theoretical Hopf algebras, again weakly group-theoretical?
\end{problem}

Affirmative answers are known under some specific  hypothesis \cite{ENO}.

\section{Composition series and length}\label{section:level-length} 

\subsection{Basic facts}\label{subsec:series-length} 
Every finite-dimensional Hopf algebra can be constructed from simple ones by succesive extensions. More precisely, we propose the following definition.

\begin{definition}\label{def:length} Let $H$ be a finite-dimensional Hopf algebra. A \emph{composition series} $\Hg$ of $H$
is a sequence of simple Hopf algebras $\Hg_1, \dots, \Hg_n$ obtained recursively as follows.

\begin{itemize}
	\item If $H$ is simple, then we let $n = 1$ and $\Hg_1 = H$.
	
	\item If $H$ is not simple, then there are $A \lhd H$, $A \neq \ku, H$, and composition series $\Ag_1, \dots, \Ag_m$, $\Bg_1, \dots, \Bg_l$, 
	of $A$ and $B = H/\hspace*{-.1cm}/A$
respectively	such that $n = m +l$ and 
\begin{align*}
\Hg_i &=	\Ag_i,& &1\le i \le m;& \Hg_i &=	\Bg_{i-m},& &m+1\le i \le m + l.
\end{align*}
\end{itemize}
\end{definition}
The simple Hopf algebras $\Hg_1$, \dots, $\Hg_n$ are the \emph{factors} of the series $\Hg$ and $n$ is its \emph{length}.
Clearly every finite-dimensional Hopf algebra admits at least one composition series.

\begin{theorem}{\cite[Theorem 1.2]{N3-JH}}\label{thm:jordan-holder} \emph{(Jordan-H\"older theorem for finite-dimen\-sional Hopf algebras).}
Let $\Hg_1, \dots, \Hg_n$ and $\Hg'_1, \dots, \Hg'_m$ be two composition series of a finite-dimensional Hopf algebra $H$.
Then there exists a bijection $\nu: \{1, \dots,n\} \to \{1, \dots,m\}$ such that $\Hg_i \simeq \Hg'_{\nu(i)}$ as Hopf algebras.
\qed
\end{theorem}

Thus, we define the \textit{length} and the \emph{factors} of $H$ as the length and the factors (up to permutation) 
of any composition series. For instance,  $H$ is of length $1$ means that it is simple; of length $2$,
that it is an extension of $K$ by $T$ where $K$ and $T$ are simple.
Also, $H$ is of length $3$ means that it fits into an extension $K\hookrightarrow H \twoheadrightarrow T$ where 
either $K$ is simple and $T$ is of length 2, or else $T$ is simple and $K$ is of length 2; but these two situations do not need to hold simultaneously,
see \cite{N3-JH}. See \cite[Section 5]{N3-JH} for the comparison of the notion of decomposition series with those of upper or lower series in \cite{MW}.

\subsubsection{Length 1}\label{subsubsec:length-1} 

Theorem \ref{thm:jordan-holder} supports the quest of simple semisimple Hopf algebras as a fundamental step
in the classification of semisimple Hopf algebras. Group algebras of simple groups, their twistings and duals are all simple \cite{Nik}, but there
are twistings of group algebras of solvable groups that are simple as Hopf algebras \cite{GN}.
However,  a twisting of the group algebra of a nilpotent group is never simple.
To our understanding, all known examples of simple semisimple Hopf algebras are twistings of group algebras or their duals.
It would be decisive either to prove that these are all, or else to find essentially new examples. We propose a working definition in this direction.

\begin{definition}\label{def:inferable} Let $H, K$ be finite-dimensional Hopf algebras. We say that $H$ is reachable from $K$ if it can be obtained from $K$
by a finite number of operations that are either duality or twisting. For instance, a cocycle-twist $K_\sigma$ is reachable from $K$:
$K \rightsquigarrow K^* \rightsquigarrow (K^*)^J \rightsquigarrow ((K^*)^J)^* = K_\sigma$, where $J$ is $\sigma$ up to identification.
Thus $H$ is reachable from $K$ if and only if either $H$ or $H^*$ is obtained from $K$ by applying succesively twists and cocycle deformations.
Clearly, being reachable is an equivalence relation.
\end{definition}

\begin{problem}\label{pbm:reachable-simple} Is every simple semisimple Hopf algebra reachable from a group algebra?
\end{problem}
In particular, we ignore the answer to:

\begin{problem}\label{pbm:reachable-simple-twist} Is a cocycle-twist of a triangular Hopf algebra again triangular?
Is there any simple Hopf algebra of the form $(\ku G^J)_{\sigma}$ but not triangular?
\end{problem}

By \cite{GN}, there are  simple  Hopf algebras  reachable from the group algebra of a super-solvable group. Also, we ask:

\begin{problem}\label{pbm:reachable-simple-classification} Classify all simple  Hopf algebras 
reachable from a group algebra.
\end{problem}

If $H$ is reachable from $K$, then they are Morita-equivalent, but it is unlikely that the converse is true.
However, we ask:

\begin{problem}\label{pbm:reachable-morita} Let $H$ and $K$ be simple semisimple Hopf algebras.
If $H$ and $K$ are Morita-equivalent, is then $H$ reachable from $K$?
\end{problem}

Perhaps the most natural, and ambitious, question is the following:

\begin{problem}\label{pbm:simple-morita} Is every \emph{simple} semisimple Hopf algebra Morita-equivalent to a group algebra?
\end{problem}

The answer to Question \ref{pbm:simple-morita} is negative without the simplicity hypothesis; indeed
Hopf algebras Morita-equivalent to a group algebra are group-theoretical. Question \ref{pbm:simple-morita} can be rephrased 
as follows-- see also Question \ref{pbm:wgt-ext}:

\begin{problem}\label{pbm:simple-morita-rephrased} 
Can every  semisimple Hopf algebra be obtained as an extension of group-theoretical ones?
\end{problem}

Finally, let us consider a Hopf algebra $H$ together with a twist $J$ and a cocycle $\sigma$, and let $H^J_\sigma$ be the vector space
$H$ with multiplication $m_{\sigma}$ and comultiplication $\Delta^J$. (This is not the same as $(H^J)_\sigma$ because being a cocycle
with respect to $\Delta$ is not the same as being a cocycle with respect to $\Delta^J$). A straightforward computation brings the 
conditions needed for $H^J_\sigma$ to be a Hopf algebra; we call this a simultaneous twist.

\begin{problem}\label{pbm:simultaneous-twist} 
Find non-trivial examples of simultaneous twists.
\end{problem}

\subsubsection{Length 2}\label{subsubsec:length-2} 
An abelian extension $\ku^{G}{\,}\leftexp{\tau}{\bowtie}_{\sigma}\ku F$ has length 2 if $G$ and $F$
are simple, but the determination of all semisimple Hopf algebras of length 2 (with known factors) is far from being clear.

\begin{problem}\label{pbm:reachable-length-2} Let $G$ and $F$ be finite non-abelian simple groups.
Find extensions of the forms
\begin{align*}
&\text{(I)}\qquad  \ku^{G} \hookrightarrow A\twoheadrightarrow \ku^{F},& &\text{(II)}\qquad \ku G\hookrightarrow A \twoheadrightarrow \ku^{F},& 
&\text{(III)}\qquad \ku G\hookrightarrow A\twoheadrightarrow \ku F,
\end{align*}
that  can not be presented as abelian extensions (in particular, they are non-trivial).
By duality, solutions to (I) give solutions to (III).
\end{problem}

\begin{example}\cite[Proposition 4.10]{N4}
 Let $G$ be a finite group and $A$ be a Hopf algebra that fits into an exact sequence 
 $\ku G\hookrightarrow A\twoheadrightarrow \ku\mathbb{Z}_2$. Then $A$ fits into an abelian exact sequence.
\end{example}

We discuss more on length 2 in the next Sections.

\section{Examples}\label{section:examples} 
Here we discuss some particular cases of extensions, see e. g. \cite[1.1, 1.2]{AN}, and specific sources of examples.
In this Section, $K$, $R$ and $T$ be Hopf algebras, while $F$, $G$, $\Gamma$, $L$ and $N$ are finite groups. 

\subsection{Smash coproduct}\label{subsec:smash-coproduct} 
First we consider the smash coproduct $R\rtimes K$; the input is a left coaction  $\rho:R\rightarrow K\otimes R$  such that
\begin{enumerate}
\item\label{1} $\Delta_R:R\rightarrow R\otimes R$, $\varepsilon:R\rightarrow\ku$, are  $K$-comodules maps;
\item\label{2} $m_R:R\otimes R\rightarrow R$, $u_R:\ku\rightarrow R$, are  $K$-comodules maps;
\item\label{3} $r_{(-1)}k\otimes r\_0=kr_{(-1)}\otimes r\_0$, for all $r\in R$, $k\in K$. 
\end{enumerate} 

Then the smash coproduct Hopf algebra $R\rtimes K$ is the tensor product algebra $R\otimes K$, where $r\#k := r\otimes k$,
with the comultiplication and antipode 
\begin{align*}
\Delta(r\# k)=r\_1\#(r\_2)_{(-1)}k\_1\otimes (r\_2)\_0\# k\_2,\,  \Ss(r\# k) = \Ss(r\_0)\#\Ss(r_{(-1)}k), 
\end{align*}
for all $r\in R$, $k\in K$. The Hopf algebra $R\rtimes K$ is an extension 
\begin{equation}\label{exact:level 2-General}
\xymatrix{ K\ar@{^{(}->}[r]^-\iota & R\rtimes K\ar@{->>}[r]^-\pi & R,}
\end{equation}
where $\iota$ and $\pi$ are  $\iota(k)=1\# k$ and $\pi(r\# k)=\varepsilon(k)r$, for  $r\in R$, $k\in K$.

\begin{remark}\label{rem:smash-coproduct-categorical}
 \cite[1.1.5]{AN} The hypothesis \eqref{1}, \eqref{2} and \eqref{3} mean that $R$ is a Hopf algebra in the category $\ydk$,
 with coaction $\rho$ and trivial action. Also, $R\rtimes K$ coincides with the bosonization $R\# K$.
\end{remark}

\begin{remark}\label{rem:smash-coproduct-dual} Assume that $\dim R < \infty$.
If $R$ is a left comodule that satisfies \eqref{1}, \eqref{2} and \eqref{3}, then $R^*$ also does, with coaction $\delta: R^* \to K\otimes R^*$ 
in the form $\langle\alpha,\Ss(r_{(-1)})\rangle \langle f,r_{(0)}\rangle = \langle\alpha,f_{(-1)}\rangle\langle f_{(0)}, r\rangle$,
for all $\alpha \in K^*$, $f \in R^*$, $r \in R$.
\end{remark}

\begin{remark}\label{rem:smash-coproduct-group}
 A  coaction $\rho:R\rightarrow\ku^G\otimes R$ satisfying \eqref{1}, \eqref{2} and \eqref{3}, 
is equivalent to a  morphism $\theta:G\rightarrow\Aut_{\Hopf} R$, 
by $\rho(r)=\sum_{\gamma\in G}\delta_\gamma\otimes\theta(\gamma^{-1})(r)$.
\end{remark}

\begin{remark}\label{rem:4}
 A left coaction $\rho:R\rightarrow\ku G\otimes R$ satisfying \eqref{1}, \eqref{2} and \eqref{3}, 
is equivalent to a $G$-grading of algebras $R=\bigoplus_{g\in G}R_g$, such that $\varepsilon$ and
$\Delta$ are homogeneous,  and $\sop(R)\subseteq Z(G)$. In particular, $\varepsilon(R_g)=0$ if $g\neq e$.
\end{remark}

\subsubsection{Extensions of group algebras} Here we assume that $R$ and $K$ are group algebras.
Let $\varphi:\ku F \to \ku G$ be an algebra isomorphism and $\psi: G \to L$ a group homomorphism.
Set $a_g= \varphi^{-1}(g)$, $g\in G$.
Then $\ku F$ has a $G$-grading of algebras $\ku F=\oplus_{g\in G}(\ku F)_g$, where $(\ku F)_g=\ku a_g$, and  an $L$-grading of algebras
$\ku F=\bigoplus_{l\in L}(\ku F)_l$, where $(\ku F)_l=\bigoplus_{g\in G:\psi(g)=l}\ku a_g$.
Let $\rho_G$, $\rho_L$  be the associated coactions and  $\rightharpoonup$  the trivial action of $\ku F$ on either 
$\ku G$ or $\ku L$.

\medbreak
$\bullet$ Let $(\rightharpoonup, \rho_L, \sigma, \tau)$ be a compatible datum. If $\psi$ is surjective, then 
the extension $\xymatrix{ \ku L\ar@{^{(}->}[r]^-\iota & \ku L {\,}\leftexp{\tau}\#_{\sigma} \ku F\ar@{->>}[r]^-\pi & \ku F}$ is abelian. 
Hence 
the extension $\xymatrix{ \ku G\ar@{^{(}->}[r]^-\iota & \ku G {\,}\leftexp{\tau}\#_{\sigma} \ku F\ar@{->>}[r]^-\pi & \ku F}$ is abelian
for any  $(\rightharpoonup, \rho_G, \sigma, \tau)$ compatible.

\pf From  \cite[(3.1.10)]{A-canad}, $(1\otimes l)\rho(a_g)=\rho(a_g)(1\otimes l)$, $\forall g\in G$, $l\in L$. Thus 
 $\sop(\ku F)=\Imm\psi =L\subseteq Z(L)$, and $L$ is abelian. \epf


$\bullet$ 
If $\varphi$ comes from an isomorphism between $F$ and $G$, then $\rho_L$ is trivial; hence $\ku F\rtimes \ku L\simeq \ku F\otimes \ku L$. 
Indeed, if $g \in G$, then $a_g\in F$, so $\varepsilon(a_g)=1$ and $a_g \in (\ku F)_{e}$, by Remark \ref{rem:4}. 
Therefore $\ku F = (\ku F)_{e}$. 
Thus we have to consider isomorphisms of group algebras not arising from group homomorphisms.
These have been intensively studied.

\medbreak
$\bullet$ In particular, we may assume that $F$, $G$ are abelian groups of the same order $n$.
Set $e_{\chi}=\frac{1}{n}\displaystyle\sum_{g\in G}\chi(g^{-1})g$, for $\chi\in \widehat{G}$; these are the primitive idempotents of $\ku G$
and 
$g=\displaystyle\sum_{\chi \in \widehat{G}} \chi(g) e_{\chi}$. Thus every isomorphism of algebras $\varphi:\ku F \to \ku G$
is determined by a bijection $\pi:\widehat{G} \to \widehat{F}$; precisely, $\varphi^{-1}: \ku G\to \ku F$ is given by 
\begin{align*}
\varphi^{-1}\left(g\right)=\dfrac{1}{n}
\displaystyle\sum_{\chi \in \widehat{G}} \sum_{x\in F} \chi(g) \pi(\chi)(x^{-1}) \, x.
\end{align*}

Then  $\rho_L$ is trivial and $\ku F\rtimes \ku L\simeq \ku F\otimes \ku L$.

\pf For each $g\in G$, \begin{align*}
\varepsilon(a_g) = \frac{1}{n} \sum_{\chi \in \widehat{G}} \sum_{x\in F} \chi(g) \pi(\chi)(x^{-1})
 =  \frac{1}{n} \sum_{\chi \in \widehat{G}} \chi(g)\, n \, \delta_{\pi(\chi), \varepsilon} = \pi^{-1}(\varepsilon) (g)\neq 0.
\end{align*} 
Thus  $a_g\in (\ku F)_{e}$ and $\rho$ is trivial.\epf

\subsection{Smash product}\label{subsec:smash-product} 

Here the input is a left action of $K$ on $T$  such that
\begin{enumerate}\setcounter{enumi}{3}
\item\label{1p} $\Delta_T:T\rightarrow T\otimes T$, $\varepsilon:T\rightarrow\ku$, are  $K$-modules maps;
\item\label{2p} $m_T:T\otimes T\rightarrow T$, $u_T:\ku\rightarrow T$, are  $K$-modules maps;
\item\label{3p} $k\_{1}\otimes k\_{2}\cdot t = k\_{2}\otimes k\_{1}\cdot t$, for all $t\in T$, $k\in K$. 
\end{enumerate} 

That is, $T$ is a Hopf algebra in $\ydk$, with the given action and trivial coaction. Then the smash product $T\# K$ is the bosonization,
i. e. the tensor product coalgebra $T\otimes K$ with the multiplication and antipode 
\begin{align*}
(t\# k)(u\# l) &= t(k\_1\cdot u) \# k\_2l,&  \Ss(t\# k) &= \Ss(k\_2)\cdot\Ss(t) \#\Ss(k\_{1}), 
\end{align*}
for all $t\in T$, $k\in K$; here $t\#k$ denotes again $t\otimes k$. The Hopf algebra $T\# K$ fits into an exact sequence of Hopf algebras 
(with obvious maps $\iota$ and $\pi$)
\begin{equation}\label{exact:level 2-General-p}
\xymatrix{ T\ar@{^{(}->}[r]^-\iota & T\# K\ar@{->>}[r]^-\pi & K.}
\end{equation}

\begin{remark}\label{rem:smash-product-dual} Assume that $\dim K < \infty$. Then
$R$ is a left $K$-comodule that satisfies \eqref{1}, \eqref{2} and \eqref{3}, if and only if $R$ a left $K^*$-module 
that satisfies \eqref{1p}, \eqref{2p} and \eqref{3p}, 
with action  $\alpha\cdot r = \langle\alpha,\Ss(r_{(-1)})\rangle r_{(0)}$, $\alpha \in K^*$, $r \in R$. Hence, assuming also that $\dim R < \infty$
and combining with Remark \ref{rem:smash-coproduct-dual}, we have
\begin{equation}\label{rem:smashcoprod-dual}
 (R\rtimes K)^* \simeq R^*\# K^*.
\end{equation}
\end{remark}

Clearly, an action of $\ku F$ on $T$ satisfying \eqref{1p}, \eqref{2p} and \eqref{3p},  is equivalent to
a morphism $\theta:F\rightarrow\Aut_{\Hopf} T$. We describe all Hopf algebra sections in
\begin{align*}
 \xymatrix@C-10pt{ T\ar@{^{(}->}[r]^-\iota & T\# \ku F\ar@{->>}[r]^-\pi & \ku F}.
\end{align*}

\begin{lemma}\label{rem:sections-smash-product}
If $\varphi\in \Hom(F, G(T))$ satisfies
\begin{align}\label{eq:section-smash}
g\cdot t &= \varphi(g) t\varphi(g^{-1}),& &g\in F, t\in T,
\end{align}
then $s_\varphi: \ku F \to  T\# \ku F$, $s_\varphi(g) = \varphi(g^{-1})\# g$, $g\in F$ is a Hopf algebra section of
$\pi$ and $K_\varphi: = \Imm s_\varphi \lhd T\# \ku F$. Moreover, any Hopf algebra section $s$
such that $\Imm s \lhd T\# \ku F$ is like this. If $\pi$ admits a Hopf algebra section with normal image, 
then $T\# \ku F\simeq T\otimes \ku F$ as Hopf algebras.
\end{lemma}

\pf If $\varphi:F\to G(T)$ is a group homomorphism satisfying \eqref{eq:section-smash}, 
then a straightforward verification shows that $s_\varphi(g)$ is a Hopf algebra section of
$\pi$ and $K_\varphi$ is normal in $T\# \ku F$.
Conversely, let $s$ be a Hopf algebra section of $\pi$ such that $K = \Imm s$ is normal in $T\# \ku F$.
Given $g\in F$, write $s (g)=\sum_{\gamma\in F} d_\gamma(g)\#\gamma$. 
Since $s$ is a coalgebra map, $\varepsilon(d_\gamma(g))=\delta_{g,\gamma}$ and $\Delta(d_\gamma(g))=d_\gamma(g)\otimes d_\gamma(g)$, for all 
$\gamma\in F$. Therefore, $s(g)=d_g(g)\# g$. Write simply $d(g)=d_g(g) \in G(T)$, for $g\in F$.
Being $s$ an algebra map, $d(gh)=d(g)(g\cdot d(h))$ and $d(1)=1$, for all $g$, $h\in F$. Since $K\lhd T\#\ku F$, 
$t_{(1)} s(g) \Ss(t_{(2)})\in K$ and $(1\#\gamma)s(g) (1\#\gamma^{-1})\in K$, for all $g$, $\gamma\in F$, $t\in T$. 
Therefore, \begin{align*}
\varepsilon(t) d(g) &=t_{(1)}d(g)(g\cdot\Ss(t_{(2)})),&   \gamma\cdot d(g)& =d(\gamma g\gamma^{-1}),& \forall t &\in T, g, \gamma \in F.
\end{align*} 
Now if $g\in F$ and $t\in T$, then
\begin{align*}
d(g)^{-1}\Ss(t)d(g)& =  d(g)^{-1} \Ss(t_{(1)}) \varepsilon(t_{(2)})d(g)\\
& = d(g)^{-1} \Ss(t_{(1)}) t_{(2)}d(g)(g\cdot\Ss(t_{(3)})) =  g\cdot \Ss(t);
\end{align*} 
hence 
$g\cdot f = d(g)^{-1}f d(g)$. Let $\varphi:F\to G(T)$, $\varphi(g)=d(g)^{-1}$ for $g\in F$. Since 
\begin{align*}
d(gh) &= d(g)(g\cdot d(h))=d(g)d(g)^{-1}d(h)d(g)=d(h)d(g),&  &\forall g, h \in F, 
\end{align*}
$\varphi$ is a group homomorphism and clearly \eqref{eq:section-smash} holds.  
Finally, observe that
\begin{align*}
s(g)(t\# 1) &= \varphi(g)^{-1}(g\cdot t)\# g=t\varphi(g)^{-1}\# g=(t\# 1)s(g),& g\in F, t\in T.
\end{align*}
Let $\psi: T\#\ku F \to T\otimes K$ defined by $\psi(t\# g)=t\varphi(g)\otimes s(g)$, $t\in T$, $g\in F$. 
Then $\psi$ is an isomorphism of Hopf algebras with inverse $\psi^{-1}:T\otimes K \to T\#\ku F$, $\psi^{-1}(t\otimes s(g))
=t\varphi(g)^{-1}\# g$, $t\in T$, $g\in F$, and the last claim follows. \epf

\subsubsection{A context}
We shall consider the following context:
\begin{itemize}
 \item $\Gamma$ is a finite simple group,
 \item $R$ is a simple semisimple Hopf algebra and
 \item  $\theta:\Gamma\rightarrow\Aut_{\Hopf} R$ is
a group homomorphism.
\end{itemize}
Then the Hopf algebras $R\rtimes \ku^\Gamma$, $R\# \ku\Gamma$, $R^*\rtimes \ku^\Gamma$, $R^*\rtimes \ku\Gamma$ are of length $2$.

\begin{lemma}\label{lemma:length2-non-abelian}
If $B\lhd R\rtimes \ku^\Gamma$ and $B \neq \ku, R\rtimes \ku^\Gamma$, then either $B = \ku^\Gamma$ or
else there exists $\varphi:\Gamma\to G(R^*)$ such that \eqref{eq:section-smash} holds and $B \simeq R$. 
\end{lemma}

\pf Set $\xymatrix@C-10pt{\ku^\Gamma\ar@{^{(}->}[r]^-\iota & R\rtimes \ku^\Gamma\ar@{->>}[r]^-\pi & R}$. 
Since $R$ is simple and $\pi(B)\lhd R$, $\pi(B)=\ku$ or $\pi(B)=R$.  Since $\Gamma$ is simple and
 $B^{\co \pi|_{B}} = B\cap (R\rtimes \ku^\Gamma)^{\co\, \pi} = B\cap \ku^\Gamma\lhd \ku^\Gamma$,  
 $B^{\co \pi|_{B}}=\ku$ or $\ku^\Gamma$. 
But $\xymatrix@C-10pt{ B^{\co \pi|_{B}}\ar@{^{(}->}[r]& B\ar@{->>}[r] & \pi(B),}$ is exact,
whence $B = \ku^\Gamma$ or $ \pi|_{B}: B \to R$ is an isomorphism of Hopf algebras. In the last case,
let $j:B \to  R\rtimes \ku^\Gamma$ be the inclusion, $p: R^*\# \ku\Gamma \to B^*$, $p = j^*$, and $K := (R^*\# \ku\Gamma)^{\co p}\lhd R^*\# \ku\Gamma$. 
Then $\iota^*|_{K}: K\to \ku\Gamma$ is an isomorphism  of Hopf algebras and Lemma \ref{rem:sections-smash-product} applies to $s = (\iota^*|_{K})^{-1}$. 
\epf

\subsection{Examples where $R$ is a twisting of a finite group}
Let $J$ be a twist of the finite group algebra $\ku N$  corresponding to a pair $(S,\omega)$.
Let $\theta :G\rightarrow\Aut_{\Hopf}(\ku N)^J$ be a morphism and 
$ H := (\ku N)^J\rtimes \ku^G$.
Note that $\wp: H\to \ku^G$, $\wp(r\# k) = \varepsilon (r)k$,
is a Hopf algebra retraction of $\iota: \ku^G \to H$. 

\begin{lemma}\label{prop:level2} Assume that 
 \begin{enumerate}\renewcommand{\theenumi}{\roman{enumi}}
\renewcommand{\labelenumi}{(\theenumi)}
\item\label{it:gamma-nonab} $G$ is non-abelian, 
\item\label{it:hyp-M} there exists $M \leqslant N^G$ non-abelian such that $[M, S] = 1$.
\end{enumerate}
Then $H$ is neither triangular nor cotriangular.
\end{lemma}
\pf
By \eqref{it:gamma-nonab},  $\ku^G$, and \emph{a fortiori} $H$, is not quasitriangular (use $\wp$). 
By \eqref{it:hyp-M}, $\ku M \leq H$; as $M$ is non-abelian,  $\ku M$ is not coquasitriangular, \emph{ditto} $H$. 
\epf

\begin{lemma}\label{comm-cocomm-biprod}
\begin{enumerate}\renewcommand{\theenumi}{\roman{enumi}}
\renewcommand{\labelenumi}{(\theenumi)}
\item\label{it:com-Biprod} $H$ is commutative if and only if $N$ is abelian. 
\item\label{it:cocom-Biprod} $H$ is cocommutative if and only if $(\ku N)^J$ is cocommutative, 
$G$ is abelian and $\rho$ is trivial.
\end{enumerate}
\end{lemma}

See Lemma \ref{lemma:cocomm Twisted group algebra} for the cocommutativity of $(\ku N)^J$.

\pf By Remark \ref{rem:smash-coproduct-categorical} and \cite[Proposition 1]{R}. 
\epf

We now show that under some assumptions, $H$ is not an abelian extension.

\begin{prop}\label{prop-examples} Assume that 
 \begin{enumerate}\renewcommand{\theenumi}{\roman{enumi}}
\renewcommand{\labelenumi}{(\theenumi)}
\item\label{it:not-abext-1} $G$ is a simple non-abelian group;
\item\label{it:not-abext-1bis} $(\ku N)^J$ is a simple Hopf algebra; 
\item\label{it:not-abext-2} if $C \lhd N$ is abelian,  then $C = \{e\}$ (in particular $N$ is non-abelian),
\item\label{it:not-abext-3} $J$ is a non-trivial twist.
\end{enumerate}
Then $H = (\ku N)^J\rtimes \ku^G$ has length 2 and is not an abelian extension.
\end{prop}

\pf $H$ is neither commutative nor cocommutative by Lemma \ref{comm-cocomm-biprod}.
Let $B \lhd H$. If $B = \ku^G$, then $H/\hspace*{-.1cm}/ B \simeq (\ku N)^J$ is not cocommutative by Lemma \ref{lemma:cocomm Twisted group algebra}
and hypotheses \eqref{it:not-abext-2} and \eqref{it:not-abext-3}.
If $B \simeq (\ku N)^J$, then it is not commutative (or alternatively $H/\hspace*{-.1cm}/ B \simeq  \ku^G$ is not cocommutative).
By Lemma \ref{lemma:length2-non-abelian}, these are the possible non-trivial normal Hopf subalgebras of $H$. 
Hence $(\ku N)^J\rtimes \ku^G$ is not an abelian extension.
\epf

\subsection{The basic construction \cite[2.1.5]{AN}}\label{section:basic2} We consider  the following data

\noindent $\bullet$ group homomorphisms $\xymatrix@C-8pt{
\Gamma\ar[r]^-\theta & \Aut_{\text{Hopf}}R &  G\ar[l]_-{\quad\mu}}$ such that $[\mu(G), \theta (\Gamma)] = 1$,

\noindent $\bullet$ a 2-cocycle $[\sigma] \in H^2 (G, \widehat \Gamma)$.

\smallbreak  
Let $A := R \rtimes \ku^{\Gamma} \#_{\sigma}\ku G$ be a Hopf algebra with underlying vector space $R \otimes \ku^{\Gamma} \otimes \ku G$ (where $r\# f \# g := r\otimes f \otimes g$)
with product, coproduct and antipode given by 
$$\aligned (r \# f \# g) (s \# f' \# g') & = r \mu (g) (s) \# f f' \sigma (g, g') \#g g', \\ 
\Delta (r \# \delta_{\gamma} \# g) & = \sum_{uv = \gamma} r\_1 \#
\delta_u \# g \otimes  \theta (u^{-1}) (r\_2) \# \delta_v \# g, \\
\Ss (r \# \delta_{\gamma} \# g) &=  \mu (g^{-1}) \circ  \theta (\gamma^{-1}) (\Ss(r)) \# \sigma (g^{-1} , g)^{-1}  \delta_{\gamma^{-1} } \#
g^{-1}, 
\endaligned $$ 
 for all $r, s \in R$, $f, f' \in \ku^{\Gamma}$, $\gamma, u, v \in \Gamma$, $g,
g' \in G$. 
It fits into an exact sequence
\begin{equation}\label{exact-A-general}
\xymatrix{
R\rtimes \ku^\Gamma \ar@{^{(}->}[r]^-\iota & R\rtimes \ku^\Gamma\#_{\sigma } \ku G\ar@{->>}[r]^-\pi & \ku G.}
\end{equation}

If $[\sigma ] = [\sigma']$  in $H^2 (G, \widehat{\Gamma})$ then $R \rtimes \ku^{\Gamma} \#_{\sigma}\ku G
\simeq R \rtimes \ku^{\Gamma} \#_{\sigma'} \ku G$.

\begin{remark} If $\Gamma$, $G$ are finite simple groups (we may assume that $\sigma$ is trivial) and $R$ is a simple semisimple Hopf algebra, 
then the Hopf algebra $R\rtimes \ku^\Gamma\# \ku G$ is of length 3.
\end{remark}

\subsubsection{Examples where $R$ a twisting of a finite group} Here we consider a basic construction
$A =(\ku N)^J\rtimes\ku^\Gamma\#\ku G$, where
 $J$ is a twist of $\ku N$  corresponding to a pair $(S,\omega)$, and  
$\theta:\Gamma\rightarrow\Aut_{\Hopf}(\ku N)^J$ and $\mu : G\rightarrow\Aut_{\Hopf}(\ku N)^J$ are group homomorphisms 
such that $[\mu(G), \theta(\Gamma)] = 1$.

\begin{prop}\label{prop-examples-2-level} Assume that:
\begin{itemize}
\item[(i)] $\Gamma$ is simple non-abelian;
\item[(ii)] $(\ku N)^J$ is a simple Hopf algebra;
\item[(iii)] if $C \lhd N$ is abelian,  then $C = \{e\}$ (in particular $N$ is non-abelian);
\item[(iv)] $J$ is a non-trivial twist.
\end{itemize} 
Then $A$ is not an abelian extension.
\end{prop}
\pf By Lemma \ref{sub-abelian-is-abelian}, since $(\ku N)^J\rtimes\ku^\Gamma$ is not so by Proposition \ref{prop-examples}.
\epf

\begin{lemma} Assume that 
 \begin{enumerate}\renewcommand{\theenumi}{\roman{enumi}}
\renewcommand{\labelenumi}{(\theenumi)}
\item $\Gamma$ is simple non-abelian, 
\item $G$ is simple,
\item $R$ is a simple semisimple Hopf algebra,
\item there exists $M \leqslant N^\Gamma$ non-abelian such that $[M, S] = 1$ and
\item\label{it:mu} $\mu(g)\in \Aut(N)$ and $(\mu(g)\otimes\mu(g))(J)=J$.
\end{enumerate}
Then $A$ is neither triangular nor cotriangular, but it fits in the exact sequence 
$\xymatrix{\ku^\Gamma\ar@{^{(}->}[r]^-j & A\ar@{->>}[r]^-p & (\ku N)^J\#\ku G}$, where $(\ku N)^J\#\ku G$ is triangular.
\end{lemma}
\pf By Lemma \ref{prop:level2}, $R\rtimes \ku^\Gamma\leqslant A$ is not coquasitriangular; then neither is $A$. 
Note that $\eta : A \to \ku^\Gamma$, $\eta(r\# f\# g)=\varepsilon(r)f$, is an epimorphism of Hopf algebras. Since $\Gamma$ is non-abelian, 
$A$ is not quasitriangular.  It is not difficult to see that the sequence  is exact. 
 By \eqref{it:mu}, $(\ku N)^J\#\ku G
\simeq (\ku N\#\ku G)^{\widetilde{J}}=(\ku(N\rtimes G))^{\widetilde{J}}$, 
where $\widetilde{J}=J_i\otimes 1 \otimes J^i\otimes 1$. 
\epf

\section{Concrete examples}\label{section:concrete examples}
In the Examples below, we consider a finite group $N$ and a twist $J$ of $\ku N$ corresponding to a pair $(S,\omega)$. Then
\begin{align}
\ce{N}{S} := \{\phi\in\Aut N:\phi|_{S}=\id_{S}\} \hookrightarrow \Aut_{\Hopf}(\ku N)^J.
\end{align}
If $Z(N) = 1$, then $\Ad: N\to \Aut N$ induces a monomorphism from the centralizer $C_N(S)$ to $\ce{N}{S}$,
and this is an isomorphism if  $\Ad: N\to \Aut N$ 
is so.
In any case, if $\Gamma \leqslant C_N(S)\geqslant G$ and $\Gamma \cap Z(N) = 1 = G \cap Z(N)$, then we denote by 
$\theta: \Gamma\rightarrow\Aut_{\Hopf}(\ku N)^J$, $\mu: G\to \Aut_{\Hopf}(\ku N)^J$ the compositions of 
the corresponding monomorphisms.

\begin{example}\label{ex:An}
Let $n, m \in \N$ such that $n\geqslant m\geqslant 9$.   Let $N = \mathbb{A}_n$,  
$\mathbb{Z}/2\times\mathbb{Z}/2\simeq S=\langle(12)(34),(13)(24)\rangle \leqslant N$, 
$0 \neq\omega\in H^2(\widehat{S},\ku^\times)$ and $J\in\ku S\otimes\ku S$ the corresponding twist. 
By \cite[4.3]{Nik}, $(\ku\mathbb{A}_n)^J$ is simple.  Let $\Gamma = \mathbb{A}_{m-4}$ (acting on $\{5,6,\cdots,m\}$); 
by hypothesis, $\Gamma$ is simple non-abelian. Now assume that $n-m \geqslant 4$; then $M=\mathbb{A}_{n-m}$ (acting on $\{m+1,m+2,\cdots,n\}$) 
is non-abelian and commutes with $\Gamma$ and $S$. By Lemma \ref{prop:level2}, 
$(\ku\mathbb{A}_n)^J\rtimes\ku^{\mathbb{A}_{m-4}}$ is neither triangular nor cotriangular. By Proposition \ref{prop-examples},  
is not an abelian extension; and is of length $2$. Now assume that $n-m \geqslant 5$; then $G=\mathbb{A}_{n-m}$ 
(acting on $\{m+1,m+2,\cdots,n\}$) is simple non-abelian and commutes with $\Gamma$. 
Therefore, $(\ku\mathbb{A}_n)^J\rtimes\ku^{\mathbb{A}_{m-4}}\# \ku \mathbb{A}_{n-m}$ is not an abelian extension but is an extension
of a triangular by a cotriangular.

\end{example}

\begin{example}\label{ex:Sn}
Let $n, m \in \N$ such that $n\geqslant m\geqslant 9$. Let $N = \mathbb{S}_n$, $\mathbb{Z}/2\times\mathbb{Z}/2\simeq S 
= \langle(12),(34)\rangle \leqslant N$, 
$0 \neq\omega\in H^2(\widehat{S},\ku^\times) \simeq \mathbb{Z}/2$ and $J\in\ku S\otimes\ku S$ the corresponding twist. 
By \cite[3.5]{GN}, $(\ku\mathbb{S}_n)^J$ is simple.
Let $\Gamma = \mathbb{A}_{m-4}$ (acting on $\{5,6,\cdots,m\}$); by hypothesis, $\Gamma$ is simple non-abelian.  
Now assume that $n-m \geqslant 3$; then $M=\mathbb{S}_{n-m}$ (acting on $\{m+1,m+2,\cdots,n\}$) is non-abelian and commutes with $\Gamma$ and $S$. 
By Lemma \ref{prop:level2}, $(\ku\mathbb{S}_n)^J\rtimes\ku^{\mathbb{A}_{m-4}}$ is neither triangular nor cotriangular; by Proposition \ref{prop-examples}, 
is not an abelian extension; and is of length $2$. Now asume that $n-m \geqslant 5$; then $G=\mathbb{A}_{n-m}$ (acting on $\{m+1,m+2,\cdots,n\}$) 
is simple non-abelian and commutes with $\Gamma$. Therefore, $(\ku\mathbb{S}_n)^J\rtimes\ku^{\mathbb{A}_{m-4}}\# \ku \mathbb{A}_{n-m}$ 
is not an abelian extension  but is an extension
of a triangular by a cotriangular.
\end{example}

\begin{example}\label{ex:PSLn} Let $\mathbb{F}_q$ be a finite field with $q$ elements. 
Let $n,r\in \N$, $s\in \N_0$ such that $n=3+r +s$, $r\geqslant 2$, $(r,q-1)=1$ and $(r,q)\neq (2,2)$. 
Let $N = PSL_n(\mathbb{F}_q)=SL_n(\mathbb{F}_q)/\{\lambda I_n:\lambda\in\mathbb{F}_q, \lambda^n=1\}$. 
Let ${\mathbb{F}_q}^\times\times {\mathbb{F}_q}^\times\simeq S =\left\{\left( \begin{matrix}
I_r &       &   &          &  \\
    & x     & 0  &     0     & \\
    &   0   & y &       0    & \\
    &   0    &  0 & (xy)^{-1} &  \\
    &        &   &          & I_s\end{matrix} \right) : x,y\in {\mathbb{F}_q}^\times\right\} \leqslant N$,   
$\omega\in H^2(\widehat{S},\ku^\times) - 0$ and $J\in\ku S\otimes\ku S$ the corresponding twist. 
Since $r\geqslant 2$, $(r,q-1)=1$ and $(r,q)\neq (2,2)$, $\Gamma=\left\{\left( \begin{matrix}
A   &         &   \\
    &   I_3    &  \\
    &          & I_s  \end{matrix}\right) :  A\in SL_r(\mathbb{F}_q)  \right\}\simeq PSL_r(\mathbb{F}_q)\leqslant N$ 
is simple non-abelian. Now assume that $s\geqslant 2$; then $M=\left\{\left( \begin{matrix}
I_r   &         &   \\
    &   I_3    &  \\
    &          & B  \end{matrix}\right) :  B\in SL_s(\mathbb{F}_q)  \right\}\simeq PSL_s(\mathbb{F}_q)\leqslant N$ is non-abelian and commutes 
    with $\Gamma$ and $S$. By Lemma \ref{prop:level2}, $(\ku PSL_n(\mathbb{F}_q))^J\rtimes\ku^{SL_r(\mathbb{F}_q)}$ is neither triangular nor cotriangular; 
    by Proposition \ref{prop-examples},  is not an abelian extension; and is of length $2$. Assume that $s\geqslant 2$, $(s,q-1)=1$ and $(s,q)\neq (2,2)$. 
Then $SL_s(\mathbb{F}_q)\simeq G=\left\{\left( \begin{matrix}
I_r   &         &   \\
      &   I_3   &   \\
       &         &  B\end{matrix}\right):  B\in SL_s(\mathbb{F}_q)  \right\}$ is simple non-abelian and commutes with $\Gamma$. Therefore, 
$(\ku PSL_n(\mathbb{F}_q))^J\rtimes\ku^{SL_r(\mathbb{F}_q)}\#\ku SL_s(\mathbb{F}_q) $ is not an abelian extension  but is an extension
of a triangular by a cotriangular.
\end{example}

\begin{remark}\label{rem:galindo} (Galindo).
Let $S\leqslant N$ abelian, $\omega\in H^2(\hat{S}, \ku^\times)$ and $J$ the corresponding twist. 
Let $G, \Gamma \leqslant C_N(S)$ such that $\Gamma \cap Z(N) = 1 = G\cap Z(N)$ and 
$\theta: \Gamma\rightarrow\Aut_{\Hopf}(\ku N)^J$ and $\mu: G\rightarrow\Aut_{\Hopf}(\ku N)^J$ as above.
Then $(\ku N)^J\rtimes \ku^\Gamma\simeq (\ku N\rtimes \ku^\Gamma)^{\widetilde{J}}$ as Hopf algebras, 
where $\widetilde{J}=J_i\otimes 1 \otimes J^i\otimes 1$. 
Since $(\ku N\rtimes \ku^\Gamma)^{\widetilde{J}}$ is group-theoretical \cite{N},
the extensions of length 2 in the Examples \ref{ex:An}, \ref{ex:Sn} and \ref{ex:PSLn} are group-theoretical. The examples of length 3 are also group-theoretical. Indeed, $(\ku N)^J\rtimes \ku^\Gamma\#\ku G\simeq \ku(N\rtimes G)^{\widetilde{J}}\rtimes \ku^\Gamma$ as Hopf algebras, where $\ku(N\rtimes G)^{\widetilde{J}}\rtimes \ku^\Gamma$ is the smash coproduct defined by $\widetilde{\theta}:\Gamma\to \Aut_{\Hopf}(\ku (N\rtimes G))^{\widetilde{J}}$, $\widetilde{\theta}(\gamma)(n,g)=(\theta(\gamma)(n),g)$, for all $\gamma\in\Gamma, n\in N, g\in G$, see Remark \ref{rem:smash-coproduct-group}. As before, $\ku(N\rtimes G)^{\widetilde{J}}\rtimes \ku^\Gamma\simeq (\ku(N\rtimes G)\rtimes \ku^\Gamma)^{\widetilde{\widetilde{J\,}}}$, where $\widetilde{\widetilde{J\,}}=J_i\otimes 1 \otimes 1\otimes J^i\otimes 1\otimes 1$, and then is group-theoretical.

\end{remark}

\begin{remark}\label{rem:davydov} (Davydov). Let $J$ be a twist of $\ku N$, $K\in \ku N\otimes\ku N$ $N$-invariant and invertible and $\phi\in \Aut N$ such that $(\phi\otimes\phi)(J)=JK$. Then $\phi\in \Aut_{\Hopf}(\ku N)^J$. In our cases, $N$ is simple and non-abelian or $N=\mathbb{S}_n$ with $n>2$, so $K\subseteq Z(N)=1$. The existence of others elements of $\Aut_{\Hopf}(\ku N)^J$ could provides new examples of extensions of length 2.
\end{remark}

\begin{problem}\label{pbm:automorphism-group-of-twist-group-algebra} What is $\Aut_{\Hopf}(\ku N)^J$?
\end{problem}


\begin{thebibliography}{ENO2}

\bibitem[A1]{A-canad} N. Andruskiewitsch, \emph{Notes on extensions of Hopf algebras}. Canad. J. Math. \textbf{48} (1996), 3--42.

\bibitem[A2]{bariloche} Andruskiewitsch, N., About finite dimensional Hopf algebras, \emph{Contemp. Math.}
\textbf{294} (2002), 1--54.

\bibitem[AD]{AD} N. Andruskiewitsch, J. Devoto, \emph{Extensions of Hopf algebras}. Algebra Anal. \textbf{7} (1995), 22--61.

\bibitem[AN]{AN} N. Andruskiewitsch, S. Natale, \emph{Examples of self-dual Hopf algebras}. J. Math. Sci. Univ. Tokyo \textbf{6} (1999), 181--215.


\bibitem[BGM]{BGM} E. Beggs, J. Gould, and S. Majid, \emph{Finite group factorizations and braiding}, J. Algebra \textbf{181} (1996), 112--151.

\bibitem[B]{B} S. Burciu, \emph{Normal Hopf subalgebras of semisimple Drinfeld doubles}. J. Pure Appl. Alg. \textbf{3:218} (2014), 540--552.

\bibitem[DT]{DT} Doi, Y., Takeuchi, M., Multiplication alteration by two-cocycles. The 
quantum version, \emph{Comm. Algebra} \textbf{22} (1994), 5715--5732. 

\bibitem[D]{D} V. Drinfeld, \emph{Quantum groups}. Proc. Int. Congr. Math., Berkeley (1987), 798--820.

\bibitem[EG1]{EG} P. Etingof, S. Gelaki, \emph{The classification of finite dimensional triangular Hopf algebras over an algebraically closed field of char 0}. Mosc. Math. J. \textbf{3} (2003), 37--43.

\bibitem[EG2]{EG2} \bysame, \emph{Some properties of finite-dimensional semisimple Hopf algebras}. Math. Res. Lett. \textbf{5} (1998), 191--197.

\bibitem[ENO1]{ENO1} P. Etingof, D. Nikshych, V. Ostrik, \emph{On fusion categories}, Ann. Math. (2) \textbf{162}, 581-642 (2005). 

\bibitem[ENO2]{ENO}  \bysame, \emph{Weakly group-theoretical and solvable fusion categories}, 
Adv. Math. \textbf{226} (2011), 176--205.

\bibitem[GNN]{GNN} S. Gelaki, D. Naidu, D. Nikshych, \emph{Centers of graded fusion categories}. Algebra Number Theory \textbf{3 (8)} (2009), 959--990.

\bibitem[GN]{GN} C. Galindo, S. Natale, \emph{Simple Hopf algebras and deformations of finite groups}. Math. Res. Lett. \textbf{14} (2007), 943--954.

\bibitem[K]{K} G.~I. Kac, \emph{Extensions of groups to ring groups}. Math. USSR Sbornik \textbf{5} (1968), 451--474.

\bibitem[LYZ]{LYZ} J.-H. Lu, M. Yan, and Y.-C. Zhu, \emph{On Hopf algebras with positive bases}, J. Algebra \textbf{237}  (2001), 421--445.

\bibitem[Ma]{M} S. Majid, \emph{Physics for algebraists: non-commutative and non-cocommutative Hopf algebras by a bicrossproduct construction}. J. Algebra \textbf{130} (1990), 17--64.

\bibitem[M]{Ma} A. Masuoka, \emph{Extensions of Hopf algebras} (lecture notes taken by Mat\'ias Gra\~na). Notas Mat. 41/99, FaMAF Uni. Nacional de C\'ordoba, 1999.


\bibitem[MW]{MW} Montgomery, S., Whiterspoon, S., \emph{Irreducible representations of
crossed products}, J. Pure Appl. Algebra \textbf{129} (1998), 315--326.

\bibitem[Mov]{Mov} M. Movshev, \emph{Twisting in group algebras of finite groups}. Func. Anal. Appl. \textbf{27} (1994), 240--244.

\bibitem[N1]{N} S. Natale, \emph{On group theoretical Hopf algebras and exact factorizations of finite groups}. J. Algebra \textbf{270} (2003), 199--211.

\bibitem[N2]{N2} \bysame, \emph{On Quasitriangular Structures in Hopf Algebras Arising from Exact Group Factorizations}. Commun. Algebra \textbf{39 (12)} (2011), 4763--4775.

\bibitem[N3]{N3-JH} \bysame, \emph{Jordan-H\"older theorem for finite dimensional Hopf algebras.} Preprint \texttt{arXiv:1407.0931},
to appear in Proc. Am. Math. Soc..

\bibitem[N4]{N4} \bysame, \emph{Hopf Algebra Extensions of Group Algebras and Tambara-Yamagami Categories}. Algebr. Represent. Th. \textbf{13} (2010), 673--691. 

\bibitem[N5]{N-pmu} \bysame, \emph{Semisimple Hopf algebras and their representations}, Publ. Mat. Uruguay \textbf{12}, 123-167 (2011).

\bibitem[Ni]{Nik} D. Nikshych, \emph{$K_0$-Rings and twisting of finite dimensional semisimple Hopf algebras}. Commun. Algebra \textbf{26 (1)} (1998), 321--342.

\bibitem[Ni2]{Nik2} \bysame, \emph{Non-group-theoretical semisimple Hopf algebras from group actions on fusion categories}. Sel. math., New ser. \textbf{14} (2008), 145--161.

\bibitem[O]{O} V. Ostrik, \emph{Module categories over the Drinfeld double of a Finite Group}, Int. Math. Res. Not.  \textbf{27}, 1507--1520 (2003). 

\bibitem[R]{R} D. Radford, \emph{The structure of Hopf algebras with a projection}. J. Algebra \textbf{92} (1985), 322--347.


\bibitem[S]{S} P. Schauenburg, \emph{Hopf bigalois extensions}. Commun. Algebra \textbf{24} (1996), 3797--3825.


\bibitem[Sch]{Sch} H.-J. Schneider, \emph{Normal basis and transitivity of crossed products for Hopf algebras}. J. Algebra \textbf{152} (1992), 289--312.

\bibitem[T]{T} M. Takeuchi, \emph{Matched pairs of groups and bismash products of Hopf algebras}. Commun. Algebra \textbf{9} (1981), 841--882.


\end{thebibliography}
\end{document}